\documentclass[11pt]{amsart}
 \usepackage{amssymb,amsfonts,amscd,amsbsy}
\usepackage[mathscr]{eucal}
\usepackage{url}

\theoremstyle{plain}
\newtheorem{thm}{Theorem}[section]

\newtheorem{prop}[thm]{Proposition}

\newtheorem{conj}[thm]{Conjecture}
\theoremstyle{definition}
\newtheorem{rem}[thm]{Remark}

\begin{document}

\title[Congruences via modular forms]
{Congruences via modular forms} 
 
\author{Robert Osburn and Brundaban Sahu}

\address{School of Mathematical Sciences, University College Dublin, Belfield, Dublin 4, Ireland}

\address{School of Mathematical Sciences, National Institute of Science Education and Research,  Bhubaneswar 751005, India}

\email{robert.osburn@ucd.ie}

\email{brundaban.sahu@niser.ac.in}

\subjclass[2000]{Primary: 11A07; Secondary: 11F11}

\date{\today}

\begin{abstract}
We prove two congruences for the coefficients of power series expansions in $t$ of modular forms where $t$ is a modular function. As a result, we settle two recent conjectures of Chan, Cooper and Sica. Additionally, we provide tables of congruences for numbers which appear in similar power series expansions and in the study of integral solutions of Ap{\'e}ry-like differential equations.
\end{abstract}
 
\maketitle

\section{Introduction}

In \cite{ccs}, Chan, Cooper and Sica investigate sequences of integers that satisfy congruence properties similar to those of the Ap{\'e}ry numbers associated with the irrationality of $\zeta(3)$. They also conjecture seven congruences and supercongruences for coefficients of power series expansions in $t$ of modular forms where $t$ is a modular function. The term supercongruences appeared in \cite{beukers1} and was the subject of the Ph.D. thesis of Coster \cite{cos}. It originally referred to families of congruences that are stronger than ones suggested by formal group theory, but now includes individual congruences (see \cite{ao}). Let   

\begin{equation*}
f(z)=\sum^\infty_{m=-\infty}\sum^\infty_{n=-\infty} q^{m^2+mn+6n^2}, \quad t_{1}=t_{1}(z)=\dfrac{\eta(z) \eta(23z)}{f(z)} 
\end{equation*}

\noindent and 

\begin{equation*}
F(z)=\sum^\infty_{m=-\infty}\sum^\infty_{n=-\infty} q^{2m^2+mn+3n^2}, \quad t_{2}=t_{2}(z)=\dfrac{\eta(z) \eta(23z)}{F(z)}
\end{equation*}

\noindent where $\eta(z)$ is the Dedekind eta-function, $q:= e^{2{\pi}iz}$ and $z \in \mathbb{H}$. Write 

\begin{equation*}
f=f(z) = \sum^\infty_{n=0} f_{n}{t_{1}^n} \quad \text{and} \quad F=F(z)=\sum^\infty_{n=0} F_{n}{t_{2}^n}. 
\end{equation*}

In \cite{ccs}, Chan, Cooper and Sica make the following

\begin{conj} \label{conjecture} If $p$ is a prime with $\left(\frac{p}{23}\right)=1$ and $n \geq 1$, then 

\begin{equation*} 
f_{np} \equiv f_{n} \pmod {p}
\end{equation*}

\noindent and

\begin{equation*} 
F_{np} \equiv F_{n} \pmod {p}.
\end{equation*}

\end{conj}

\noindent The first few terms in the sequence $\{f_{n}\}_{n\geq 0}$ are

\begin{center}
$1$, $2$, $6$, $26$, $142$, $876$, $5790$, $40020$, $285582$, $\dots$
\end{center}

\noindent while for $\{F_{n}\}_{n\geq 0}$, we have

\begin{center}
$1$, $0$, $2$, $6$, $30$, $144$, $758$, $4080$, $22702$, $128832$, $\dots$.
\end{center}

\noindent Closed forms for $f_{n}$ and $F_{n}$ were not given in \cite{ccs} and thus a combinatorial approach to Conjecture \ref{conjecture} is not yet available. The purpose of this note is to prove this conjecture via modular forms. We have the following.

\begin{thm} \label{main} If $p$ is a prime with $\left(\frac{p}{23}\right)=1$ and $n$, $r \geq 1$ are integers, then

\begin{equation} \label{conj1}
f_{np^r} \equiv f_{np^{r-1}} \pmod {p^r}
\end{equation}

\noindent and 

\begin{equation} \label{conj2}
F_{np^r} \equiv F_{np^{r-1}} \pmod {p^r}.
\end{equation}
\end{thm}

In Section 2, we recall some preliminaries on power series expansions and Eisenstein series with characters and then prove Theorem \ref{main}. In Section 3, we provide tables of congruences for numbers appearing in other power series expansions found in \cite{ccs} and in the study of integral solutions of Ap{\'e}ry-like differential equations (see \cite{avsz}, \cite{beukers3}, \cite{zagier}) and mention conjectural supercongruences. Finally, we note that two other conjectural congruences from \cite{ccs} which involve $f_{2,n}$ and $f_{3,n}$ (see Section 3) have recently been proven in \cite{ckko}. The remaining three conjectures in \cite{ccs} are still open. 

\section{Proof of Theorem \ref{main}}

We first recall a recent result of Jarvis and Verrill (see Proposition 4.2 in \cite{jv} or Proposition 3 in \cite{beukers2}). This result is quite useful as it allows one to deduce congruence properties of coefficients in a power series expansion from those of another expansion.

\begin{prop}\label{jv}
Let $t$ be a power series 
$$
t=\frac{1}{m} \sum^{\infty}_{n=1} a_nu^{n/v},
$$
convergent in a neighborhood of $u=0,$ with $m, v$ positive integers, $a_n \in \mathbb Z$ and 
$a_1=1.$ Suppose that in some neighborhood of $u=0$ we have an equality of convergent power
series given by 

\begin{equation} \label{equal}
\sum^{\infty}_{n=1} b_nt^{n-1} \, dt= \sum^{\infty}_{n=1} c_nu^{n-1} \, du,
\end{equation}

\noindent for some integers $b_n$ and $c_n,$ $n \ge 1$. Assume $p$ is a prime not dividing $m$ or $v.$ If
$$
b_{np^r} \equiv  b_{np^{r-1}} \pmod {p^r},
$$
then 
$$
c_{np^r} \equiv  c_{np^{r-1}} \pmod {p^r}.
$$
  
\end{prop}

\begin{rem} \label{stien}
J. Stienstra has kindly pointed out that one can use formal group theory (see the appendix of \cite{sb}) to extend Proposition \ref{jv} to the case where each of the sums in (\ref{equal}) starts with $n=0$. Also, since $a_1=1$, the converse of Proposition \ref{jv} is true.
\end{rem}

We now discuss the notion of Eisenstein series with characters. For further details, see Chapter $5$ of \cite{stein}. Let $M_{k}(\Gamma_{0}(N), \epsilon)$ be the space of modular forms of weight $k$ on $\Gamma_{0}(N)$ with character $\epsilon$. Suppose $\chi$ and $\psi$ are primitive Dirichlet characters with conductors $L$ and $R$, respectively. Let

\begin{equation} \label{ek}
E_{k, \chi, \psi}(q) := c_{0} + \sum_{n=1}^{\infty} \Bigl( \sum_{d \mid n} \psi(d) \chi(n/d) d^{k-1} \Bigr) q^n
\end{equation}

\noindent where 

$$
c_{0}=\left\{ \begin{array}{ll} -\dfrac{B_{k, \psi}}{2k} & \qquad \textup{if} \quad L=1, \vspace{.05in} \\
\qquad 0 & \qquad \textup{if} \quad L>1 \end{array} \right.
$$

\noindent and $B_{k, \psi}$ is the generalized Bernoulli number associated to $\psi$. If $t$ is a positive integer and $k \geq 3$ is an integer such that $\chi(-1)\psi(-1)=(-1)^k$, then $E_{k, \chi, \psi}(q^t)$ is in $M_{k}(\Gamma_{0}(RLt), \chi \psi)$. Moreover, given $N$ and $\epsilon$, the series $E_{k, \chi, \psi}(q^t)$ such that $RLt \mid N$ and $\chi \psi =\epsilon$ form a basis for the Eisenstein subspace $E_{k}(\Gamma_{0}(N), \epsilon)$ of $M_{k}(\Gamma_{0}(N), \epsilon)$. \\

\noindent {\it{Proof of Theorem \ref{main}.}} 
Let $\chi$ be the character $\left(\frac{\cdot}{23}\right)$ and $\psi$ be the trivial character $1$. 
We first note that 

$$E_{3, 1, \chi}(q) =: \sum_{n=0}^{\infty} e_n q^{n}$$

\noindent and

$$E_{3, \chi, 1}(q)=:\sum_{n=0}^{\infty} a_n q^n$$ 

\noindent form a basis for the space $E_3\Bigl(\Gamma_0(23), \left(\frac{\cdot}{23}\right)\Bigr)$. By Lemma 0.3 in \cite{verrill} and a finite computation, we have

\begin{equation*}
\displaystyle f \frac{q\frac{dt_1}{dq}}{t_1} = F \frac{q\frac{dt_2}{dq}}{t_2}= -\frac{1}{24} E_{3, 1, \chi}(q) - \frac{23}{24} E_{3, \chi, 1}(q)
\end{equation*}

\noindent and so 

\begin{equation} \label{pre}
f \frac{dt_1}{t_1} = F \frac{dt_2}{t_2}= \Biggl[-\frac{1}{24} E_{3, 1, \chi}(q) - \frac{23}{24} E_{3, \chi, 1}(q) \Biggr] \dfrac{dq}{q}.
\end{equation}

\noindent By (\ref{ek}), we have

\begin{equation} \label{co1}
a_{np^r} -  \left(\frac{p}{23}\right) a_{np^{r-1}} = \sum_{d^{\prime} \mid n}  \left(\frac{n/d^{\prime}}{23}\right) (d^{\prime} p^r)^2
\end{equation}

\noindent and

\begin{equation} \label{co2}
e_{np^r} - e_{np^{r-1}} = \sum_{d^{\prime} \mid n}  \left(\frac{d^{\prime}p^r}{23}\right) (d^{\prime} p^r)^2.
\end{equation}

\noindent Letting $u=q$ in (\ref{pre}) implies that 

\begin{equation} \label{setup}
\sum_{n=0}^{\infty} f_n t_1^{n-1} \, dt_1 = \sum_{n=0}^{\infty} F_n t_2^{n-1} \, dt_2 = \Biggl[-\frac{1}{24} E_{3, 1, \chi}(u) - \frac{23}{24} E_{3, \chi, 1}(u) \Biggr] \dfrac{du}{u}.
\end{equation}

\noindent If we take $v=1$, $m=1$, $b_n=f_n$, $F_n$, respectively, and $c_n=\displaystyle -\frac{1}{24} e_{n} - \frac{23}{24} a_{n}$, then by (\ref{co1}) and (\ref{co2}), we have 

\begin{equation} \label{check}
c_{np^r} \equiv c_{np^{r-1}} \pmod {p^r}
\end{equation} 

\noindent for primes  $p \geq 3$ such that $\left(\frac{p}{23}\right)=1$ and $r \geq 1$ and for $p=2$ and $r \geq 3$. An application of Remark \ref{stien} then implies (\ref{conj1}) and (\ref{conj2}). To verify (\ref{check}) for $p=2$ and $r=1$ or $2$, we first note that (\ref{ek}) implies $e_{n}=a_{n}$ if $\left(\frac{n}{23}\right)=1$ and  $e_{n}=-a_{n}$ if $\left(\frac{n}{23}\right)=-1$. The case where $\left(\frac{n}{23}\right)=0$ can be reduced to one of previous two cases since $a_{23n}=23^2 a_{n}$ and $e_{23n}=e_{n}$ for all $n$. The result then follows upon a routine check that (\ref{check}) holds in all of these cases.

\qed

\section{Tables}

Using the methods in Section 2, we have proven congruences of the form

\begin{equation} \label{an} 
A(np^r) \equiv A(np^{r-1}) \pmod {p^r}
\end{equation}

\noindent for all of the numbers $A(n)$ which appear in Tables 1, 2 and 3. For brevity, we only give the relevant modular function $t$ and modular forms $f(t)$ and $\displaystyle M:=f(t) \frac{q\frac{dt}{dq}}{t}$. The coefficients of $M$ in Tables 1, 2 and 3 can be computed using (\ref{ek}), Chapter 4, Section 32 in \cite{fine} (for example, see page 85, equation (32.71) for (vi)), \cite{wom} or \cite{stein}. Given positive integers $s_{1}$, $s_{2}$, $\dotsc$, $s_{k}$ and integers $r_{1}$, $r_{2}$, $\dotsc$, $r_{k}$, we write 

$$s_{1}^{r_1} s_{2}^{r_2} \dotsc s_{k}^{r_k}$$ 

\noindent for the eta-quotient 

$$\eta(s_{1} z)^{r_1}\eta(s_{2} z)^{r_2} \cdots \eta(s_{k} z)^{r_k}.$$ 

Table 1 consists of numbers $f_{i,n}$, $i=2$, $3$, $5$, $7$ and $11$ which are coefficients in the power series expansion in $t$ of the modular forms (see \cite{ccs})

\begin{equation*}
f_{2}=\sum^\infty_{m=-\infty}\sum^\infty_{n=-\infty} q^{m^2+n^2}, \quad f_{3}=\sum^\infty_{m=-\infty}\sum^\infty_{n=-\infty} q^{m^2+mn+n^2}, \quad f_{5}= \dfrac{1^5}{5^1},
\end{equation*}

\begin{equation*}
 f_{7}=\sum^\infty_{m=-\infty}\sum^\infty_{n=-\infty} q^{m^2+mn+2n^2} \quad \text{and} \quad  f_{11}=\sum^\infty_{m=-\infty}\sum^\infty_{n=-\infty} q^{m^2+mn+3n^2}.
\end{equation*}

\noindent We write $\chi_{s}:=\left(\frac{\cdot}{s}\right)$ for $s=3$, $5$, $7$, $11$ and $\chi_{-4}:=\left(\frac{-4}{\cdot}\right)$. The analogue of Theorem \ref{main} is true for primes $p$ satisfying $\chi_{-4}(p)=1$ in (i), $\chi_{3}(p)=1$ in (ii), $\chi_{7}(p)=1$ in (iv), $\chi_{11}(p)=1$ in (v). It is true for all primes in (iii). The only known closed forms are 

$$f_{2,n}=\left(\frac{8^n\left(\frac{1} {4}\right)_n}  {n!}\right)^2$$ 

\noindent and 

$$f_{3, n}=\frac{108^n\left(\frac{1} {6}\right)_n \left(\frac{1} {3}\right)_n}  {(n!)^2}.$$ 

Table 2 lists numbers which arise in Beukers' \cite{beukers3} and Zagier's \cite{zagier} study of integral solutions of second order Ap{\'e}ry-like differential equations. The choices of $t$ and the parameterizations of $f$ can be found in \cite{verrill} and \cite{zagier}. In case (ix), congruence (\ref{an}) with $A(n)$ replaced by $2^{n} A(n)$ has been proven in \cite{jv}. 

Table 3 contains numbers listed in \cite{avsz} as part of a discussion on third order Ap{\'e}ry-like differential equations. Here 

$$
L_{1}(z):=-\dfrac{7}{240} E_4(z) +\dfrac{1}{60} E_4(2z)-\dfrac{3}{80} E_4(3z) + \dfrac{21}{20} E_4(6z), 
$$

$$
L_{2}(z):= \dfrac{1}{120} E_4(z) -\dfrac{2}{15} E_4(2z) -\dfrac{3}{40} E_4(3z)  +\dfrac{6}{5} E_4(6z),
$$

\noindent and

$$
L_{3}(z):=\dfrac{1}{240} E_4(z) -\dfrac{1}{60} E_4(2z) -\dfrac{27}{80} E_4(3z) +\dfrac{27}{20} E_4(6z)
$$

\noindent where $E_{4}(z)$ is the usual weight $4$ Eisenstein series on $SL_{2}(\mathbb{Z})$. The choices of $t$ and the parameterizations of $f$ can be found in \cite{ccl}, \cite{cv} and \cite{ps}. 

Finally, we have numerically observed extensions of (\ref{an}) modulo $p^{2r}$ (subject to the above conditions for $p$ odd) in (i), (ii), (iv), (v), (vii), (viii), (ix) and (x) and modulo $p^{3r}$ for (iii), (xii) and (xiii). Here $p \geq 5$ for (xii). Coster \cite{cos} has proven an extension of (\ref{an}) modulo $p^{3r}$ for (vi) and (xi). It might of interest to see if combinatorial techniques can be applied to some of these conjectural extensions.

\begin{center}
\begin{table}
\caption{}
\begin{tabular}{|c|c|c|c|}
\multicolumn{3}{r}{} \\ \hline
$A(n)$ & $t$ & $M$ \\
\hline 
& & \\
(i) \hspace{.05in} $f_{2,n}$ & $\dfrac{2^{12}}{f_{2}^6}$ & $-4E_{3, 1, \chi_{-4}}(q) - 16E_{3, \chi_{-4}, 1}(q)$ \\
& & \\
(ii) \hspace{.05in} $f_{3,n}$ & $\dfrac{1^6 3^6}{f_{3}^6}$ & $-9E_{3, 1, \chi_{3}}(q) - 27E_{3, \chi_{3}, 1}(q)$  \\
& & \\
(iii) \hspace{.05in} $f_{5,n}$ & $\dfrac{5^6}{1^6}$ &  $E_{4, 1, \chi_{5}}(q)$ \\
& & \\
(iv) \hspace{.05in} $f_{7,n}$ & $\dfrac{1^3 7^3}{f_{7}^3}$ & $-\frac{7}{8} E_{3, 1, \chi_{7}}(q) - \frac{49}{8} E_{3, \chi_{7}, 1}(q)$ \\ 
& &  \\
(v) \hspace{.05in} $f_{11, n}$ & $\dfrac{1^2 11^2}{f_{11}^2}$ & $-\frac{1}{3} E_{3, 1, \chi_{11}}(q) - \frac{11}{3} E_{3, \chi_{11}, 1}(q)$ \\
&& \\
\hline

\end{tabular} 
\end{table}
\end{center}

\begin{center}
\begin{table} 
\caption{}
\begin{tabular}{|c|c|c|c|c|}
\multicolumn{3}{r}{} \\ \hline
$A(n)$ & $t$ & $f(t) $ & $M$ \\
\hline 
& & & \\
(vi) \quad $\displaystyle \sum^n_{k=0} \binom {n} {k}^3$ & $\dfrac{1^3 6^9}{2^3 3^9}$ & $\dfrac{2^1 3^6}{1^2 6^3}$ & $\dfrac{1^1 2^4 3^5}{6^4}$ \\
& & & \\
(vii) \quad $\displaystyle \sum^{\lfloor \frac{n}{3} \rfloor}_{k=0} (-1)^k 3^{n-3k}\binom {n} {3k} \binom {3k} {k} \binom{2k}{k}$ & $-\dfrac{9^3}{1^3}$  & $\dfrac{1^3}{3^1}$ & $\dfrac{3^9}{9^3}$ \\
& & & \\
(viii) \quad $\displaystyle \sum^n_{k=0} \binom {n} {k}^2\binom{2k}{k}$ & $\dfrac{1^4 6^8}{2^8 3^4}$ & $\dfrac{2^6 3^1}{1^3 6^2}$ &$\dfrac{1^1 2^4 3^5}{6^4}$ \\
& & & \\
(ix) \quad $\displaystyle \sum^{\lfloor \frac{n}{2} \rfloor}_{k=0} 4^{n-2k}\binom {n} {2k} \binom {2k} {k}^2$ & $\dfrac{1^4 4^2 8^4}{2^{10}}$ & $\dfrac{2^{10}}{1^4 4^4}$ & $\dfrac{2^4 4^6}{8^4}$ \\
& & &  \\
(x) \quad $\displaystyle \sum_{k=0}^{n} \sum_{l=0}^{k} (-1)^{k} 8^{n-k}\binom {n} {k} \binom {k} {l}^3$ & $\dfrac{1^5 3^1 4^5 6^2 12^1}{2^{14}}$ & $\dfrac{2^{15} 3^2 {12}^2}{1^6 4^6 6^5}$ &$\dfrac{2^{7} 6^{11}}{1^1 3^5 4^1 12^5}$ \\
& & &  \\
\hline
\end{tabular} 
\end{table}
\end{center}

\begin{center}
\begin{table} 
\caption{}
\begin{tabular}{|c|c|c|c|}
\multicolumn{3}{r}{} \\ \hline
$A(n)$ & $t$ & $f$ & $M$ \\
\hline 
& & & \\
(xi) \quad $\displaystyle \sum_{k=0}^{n} \binom{n+k}{k}^2 \binom{n}{k}^2$ & $\dfrac{1^{12} 6^{12}}{2^{12} 3^{12}}$ & $\dfrac{2^7 3^7}{1^5 6^5}$ & $L_1(z)$ \\ 
& & & \\
(xii) \quad $\displaystyle (-1)^n\sum_{k=0}^{n} \binom{n}{k}^2 \binom{2k}{k} \binom{2(n-k)}{n-k}$ & $\dfrac{2^6  6^6}{1^6 3^6}$ & $\dfrac{1^4 3^4}{2^2 6^2}$ & $L_2(z)$ \\
& & & \\
(xiii) \quad $\displaystyle  (-1)^n\sum_{k=0}^{[n/3]} (-1)^k \frac{ 3^{n-3k}(3k)!}{(k!)^3}  \binom{n}{3k}\binom{n+k}{k}$ & $\dfrac{3^4 6^4}{1^4 2^4}$ & $\dfrac{1^3  2^3}{3^1 6^1}$ & $L_3(z)$ \\
& & & \\
\hline
\end{tabular} 
\end{table}
\end{center}

\section*{Acknowledgements}
The authors would like to thank Heng Huat Chan, Shaun Cooper, Jan Stienstra and the referee for their careful reading of the paper and for many helpful comments. The second author thanks the Institut des Hautes {\'E}tudes Scientifiques for their hospitality and support during the preparation of this paper. The authors were partially supported by Science Foundation Ireland 08/RFP/MTH1081.


\begin{thebibliography}{10}

\bibitem{ao}
S. Ahlgren, K. Ono, \emph{A Gaussian hypergeometric series evaluation and Apery number congruences}, J. Reine Angew. Math. \textbf {518}, (2000), 187--212. 

\bibitem{avsz}
G. Almkvist, D. van Straten and W. Zudilin, \emph{Generalizations of Clausen's formula and algebraic transformations of Calabi-Yau differential equations}, Proc. Edinb. Math. Soc. (2), to appear.

\bibitem{beukers1}
{F. Beukers}, {\em Some congruences for the Ap{\'e}ry numbers}, {J. Number Th.} {\bf 21} (1985), no. 2, 141--155.

\bibitem{beukers2}
{F. Beukers}, {\em Another congruence for the Ap{\'e}ry numbers}, {J. Number Th.} {\bf 25} (1987), no. 2, 201--210.

\bibitem{beukers3}
F. Beukers, \emph{On B. Dwork's accessory parameter problem}, Math. Z. \textbf{241} (2002), no. 2, 425--444.

\bibitem{ccl}
H. Chan, S. Chan and Z. Liu, \emph{Domb's numbers and Ramanujan-Sato type series for $1/\pi$}, Adv. Math. \textbf{186} (2004), 396--410.

\bibitem{ccs}
H. Chan, S. Cooper and F. Sica, \emph{Congruences satisfied by Ap{\'e}ry-like numbers}, Int. J. Number Theory \textbf{6} (2010), no. 1, 89--97.

\bibitem{ckko}
H. Chan, A. Kontogeorgis, C. Krattenthaler and R. Osburn, \emph{Supercongruences satisfied by coefficients of ${_2F_1}$ hypergeometric series}, Ann. Sci. Math. Qu{\'e}bec \textbf{34} (2010), no. 1, 25--36.

\bibitem{cv}
H. Chan, H. Verrill, \emph{The Ap{\'e}ry numbers, the Almkvist-Zudilin numbers and new series for $1/\pi$}, Math. Res. Lett. \textbf{16} (2009), no. 3, 405--420.

\bibitem{cos}
M. Coster, \emph{Supercongruences}, Ph.D. thesis, Universiteit Leiden, 1988.

\bibitem{fine}
N. Fine, \emph{Basic hypergeometric series and applications}, American Mathematical Society, Providence, RI, 1988.

\bibitem{jv} 
F. Jarvis, H. Verrill, \emph{Supercongruences for the Catalan-Larcombe-French numbers}, Ramanujan J. \textbf{22} (2010), no. 2, 171--186.

\bibitem{wom}
K. Ono, \emph{The web of modularity: arithmetic of the coefficients of modular forms and $q$-series}, Amer. Math. Soc., CBMS Regional Conf. in Math., vol. 102, 2004.

\bibitem{ps}
C. Peters, J. Stienstra, \emph{A pencil of $K3$-surfaces related to Ap{\'e}ry's recurrence for $\zeta(3)$ and Fermi surfaces for potential zero}. Arithmetic of complex manifolds (Erlangen, 1988), 110--127, Lecture Notes in Math., 1399, Springer, Berlin, 1989.

\bibitem{stein}
W. Stein, \emph{Modular forms, a computational approach}, Graduate Studies in Mathematics, \textbf{79}. American Mathematical Society, Providence, RI, 2007.

\bibitem{sb}
J. Stienstra, F. Beukers, \emph{On the Picard-Fuchs equation and the formal Brauer group of certain elliptic $K3$ surfaces}, Math. Ann. \textbf{271} (1985), no. 2, 269--304.

\bibitem{verrill}
H. Verrill, {\em Some congruences related to modular forms}, available at \url{http://www.mpim-bonn.mpg.de/Research/MPIM+Preprint+Series}

\bibitem{zagier}
D. Zagier, {\em Integral solutions of Ap\'ery-like recurrence equations}, Group and Symmetries: From Neolithic Scots to John McKay, 349--366, CRM Proc. Lecture Notes, \textbf{47}, Amer. Math. Soc., Providence, RI, 2009.

\end{thebibliography}
\end{document}